\documentclass{article}

 \usepackage{color}
 
\usepackage{tikz}
\usepackage{amsmath}
\usepackage{amsfonts}
\usepackage{amssymb}
\usepackage{graphicx}
\usepackage{array}
\usepackage{fullpage}

\newcommand*\circled[1]{\tikz[baseline=(char.base)]{
            \node[shape=circle,draw,inner sep=2pt] (char) {#1};}}

\begin{document}
%
%
%
%
%

\title{\textsc{A method for identifying stability regimes using roots of a reduced-order polynomial}}
\author{Olga Trichtchenko}
\date{}



\maketitle

\begin{abstract}
For dispersive Hamiltonian partial differential equations of order $2N+1$, $N\in \mathbb{Z}$, there are two criteria to analyse to examine the stability of small-amplitude, periodic travelling wave solutions to high-frequency perturbations. The first necessary condition for instability is given via the dispersion relation. The second criterion for instability is the signature of the eigenvalues of the spectral stability problem given by the sign of the Hamiltonian. In this work, we show how to combine these two conditions for instability into a polynomial of degree $N$. If the polynomial contains no real roots, then the travelling wave solutions are stable. We present the method for deriving the polynomial and analyse its roots using Sturm's theory via an example. 
\end{abstract}

\section{Introduction}
Partial differential equations (PDEs) are used in a wide variety of applications to describe physical phenomena where this physical relevance imposes the requirement that the solutions to the PDEs are real. Moreover, if the description is of a closed system, there is usually an associated conservation of energy and the equations used are Hamiltonian. As more methodology for solving PDEs is developed \cite{OlverBook}, the natural question to ask is then how realistic are these solutions are and how likely we are to observe them in nature. Thus, analysing their stability also becomes important \cite{BJK11, hurjohnson, KM14}. The purpose of this work is to present a simplified method for stability analysis, illustrated by an explicit example. We focus on high-frequency instabilities arising from spectral analysis of a perturbation of periodic travelling waves \cite{DT15} and restrict our focus to stability of solutions of dispersive Hamiltonian equations. We show how working with the dispersion relation, we can methodically construct a parameter regime where there is only spectral stability with respect to particular perturbations and in the regions where we expect instability, we show what types of instabilities can arise. 

In recent work \cite{DT15}, a method for establishing the presence of high-frequency instabilities of travelling wave solutions for both scalar PDEs as well as for systems of equations was described. In this method, there are two important conditions to consider: 
\begin{enumerate}
\item collisions of eigenvalues of the spectral stability problem and
\item the signature of these eigenvalues.
\end{enumerate} 
Furthermore, it was shown that in order for the solutions to become unstable, the system had to admit waves travelling in different directions (bi-directional waves). In the follow-up work \cite{TDK18}, the authors showed that a different way to meet the instability criteria, was for equations to contain what is referred to as a generalised resonance. An equation contains a resonance if there is a certain set of parameters for which travelling wave solutions are predominantly composed of at least two distinct frequencies which can travel at the same speed. Physically, this implies that there are at least two different forces that can influence the travelling waves that are of the same order of magnitude. For example, if we are considering water waves, then these waves are in a \textbf{resonant regime} if surface tension and gravity are competing forces of the same order of magnitude. The result is that the travelling wave profiles contain two different prominent modes, otherwise referred to as Wilton ripples \cite{TDW16, W15}.

If we restrict ourselves to scalar, dispersive and Hamiltonian PDEs where the solution $u$ depends on one spatial and one time variable, i.e. $u = u(x,t)$ with a period $L$ and up to $2N+1$ derivatives, then it has been shown \cite{DT15} that all we need is a polynomial dispersion relation $\omega(k)$ of order $2N+1$ to describe both of the necessary conditions for instability. In \cite{KDT18}, it was shown that the two necessary conditions for instability can be collapsed into one criterion on the roots of a polynomial of order $N$ to be in an interval $I$ defined in Section \ref{sec:stabTheory}. This greatly simplifies the analysis, leading to closed-form results for stability regions of specific PDEs. 

This work presents a method for the single criteria for instability of periodic travelling wave solutions to a dispersive, Hamiltonian PDE using an example with three competing terms. The formulation and underlying theory is described in Section \ref{sec:stabTheory}. Working with the dispersion relation, we show the general methodology for the stability analysis in Section \ref{sec:methods}. Section \ref{sec:example} explicitly shows how to implement the method via an example, demonstrating how to construct the coefficients systematically and use Sturm's theory to analyse the roots of the reduced polynomial. In Section \ref{sec:results}, figures of the stability and instability regions are shown and we conclude in Section \ref{sec:conclusion}.

\section{Summary of Stability Theory}\label{sec:stabTheory}
Consider a scalar Hamiltonian PDE of the form
\begin{align}
u_t = \partial_x \frac{\delta H}{\delta u},
\end{align} 
where the function $u = u(x,t)$ describes a periodic travelling wave, with $H$ the Hamiltonian and $\frac{\delta H}{\delta u}$ a variational derivative. More specifically $u(x,t)$ is a solution of
\begin{align}
u_t = \sum_{n=1}^{N} C_{2n+1} \frac{\partial^{2n+1}u}{\partial x^{2n+1}} +  f(u,u_x,...,u_{(2N)x})_x,
\label{eq:genForm}
\end{align}
where $N$ is positive integer and $\frac{\partial^{2n+1}u}{\partial x^{2n+1}}$ are $2n+1$ (odd) derivatives up to order $2N+1$ with the nonlinearity $f$ that can depend on $u$ as well as its derivatives up to order $2N$ (denoted as $u_{(2N)x}$), keeping the overall system dispersive. For ease, we consider the equation with real coefficients $C_{2n+1}$. We obtain the dispersion relation $\omega(k)$ if we let $u(x,t) \sim e^{ikx - i\omega t}$ with $k$ a Fourier mode, and substitute into \eqref{eq:genForm} to obtain
\begin{align}
\omega(k) = \sum_{n=1}^N (-1)^{(n+1)}C_{2n+1}k^{2n+1}.
\label{eq:dispersion}
\end{align}
Furthermore, if we restrict the space of solutions $u(x,t)$ to periodic, travelling waves moving at speed $V$ such that $u(x,t)\rightarrow u^{(0)}(x-Vt)$, then we can write \eqref{eq:genForm} in the travelling frame of reference and consider the steady-state equation
\begin{align}
 Vu_x + \sum_{n=1}^{N} C_{2n+1} \frac{\partial^{2n+1}u}{\partial x^{2n+1}} + f(u,u_x,...,u_{(2N)x})_x = 0,
\label{eq:travGenForm}
\end{align}
and setting $x\rightarrow x-Vt$ from now on. Despite restricting the space of solutions to travelling waves $u^{(0)}(x)$, we can still gather information about the time dependence by perturbing about this steady-state with a small perturbation governed by $\delta$, i.e.
\begin{align}
u(x,t) & = u^{(0)}(x) + \delta \bar{u}^{(1)}(x,t) \nonumber \\
& = u^{(0)}(x) + \delta e^{\lambda t} u^{(1)}(x).
\label{eq:deltaExp}
\end{align}
We have made an assumption about the time dependence of the perturbation by introducing $\lambda \in \mathbb{C}$. Recall that $u^{(0)}(x)$ is periodic of period $L$ (for convenience, $L=2\pi$) \cite{DT15}. We allow the perturbations to be of any period, but bounded in space using the Fourier-Floquet expansion
\begin{align}
u^{(1)}(x) = e^{i\mu x}\sum_{m=-M}^{M} b_{m} e^{imx},
\label{eq:pertForm}
\end{align}
with $\mu \in \mathbb{R}$ the Floquet parameter governing the period of the perturbation and a Fourier mode $m \in \mathbb{Z}$ \cite{DK06}. We note that this perturbation can grow exponentially in time if Re$(\lambda)>0$, where $\lambda = \lambda(\mu+m)$ depends on the Fourier-Floquet modes $m$ and $\mu$. For solutions with $|u^{(0)}(x)| = O(\epsilon)$ with $\epsilon \rightarrow 0$, 
\begin{align}
\lambda(\mu+m) =  i(m+\mu)V -i\omega(m+\mu),
\end{align}
if we consider $O(\delta)$ term when substituting \eqref{eq:deltaExp} and \eqref{eq:pertForm} into \eqref{eq:genForm}, staying in the travelling frame of reference.

\begin{figure}
\begin{center}
\begin{tikzpicture}
\draw[help lines, color=gray!30, dashed] (-1.1,-2.1) grid (1.1,2.1);
\draw[->,ultra thick] (-1.0,0)--(1.0,0) node[right]{Re$\lambda$};
\draw[->,ultra thick] (0,-2)--(0,2) node[above]{Im$\lambda$};
\draw[color=blue, fill] (0,1) circle (0.1cm);
\draw[color=blue, fill] (0,1.5) circle (0.1cm);
\draw[color=blue, fill] (0,-1) circle (0.1cm);
\draw[color=blue, fill] (0,-1.5) circle (0.1cm);
\end{tikzpicture}
\begin{tikzpicture}
\draw[help lines, color=gray!30, dashed] (-1.6,-2.1) grid (1.6,2.1);
\draw[->,ultra thick] (-1.5,0)--(1.5,0) node[right]{Re$\lambda$};
\draw[->,ultra thick] (0,-2)--(0,2) node[above]{Im$\lambda$};
\draw[color=red, fill] (1,1) circle (0.1cm);
\draw[color=red, fill] (1,-1) circle (0.1cm);
\draw[color=red, fill] (-1,1) circle (0.1cm);
\draw[color=red, fill] (-1,-1) circle (0.1cm);
\end{tikzpicture}
\begin{tikzpicture}
\draw[help lines, color=gray!30, dashed] (-1.6,-2.1) grid (1.6,2.1);
\draw[->,ultra thick] (-1.5,0)--(1.5,0) node[right]{Re$\lambda$};
\draw[->,ultra thick] (0,-2)--(0,2) node[above]{Im$\lambda$};
\draw[color=red, fill] (1.2,0) circle (0.1cm);
\draw[color=red, fill] (-1.2,0) circle (0.1cm);
\end{tikzpicture}
\caption{Three different configurations of the smallest number of eigenvalues $\lambda$ of the spectral stability problem of a Hamiltonian system, showing the symmetry about the real and imaginary axes. On the left (in blue), is the stable regime. The centre and right panel are the unstable regimes (in red).\label{fig:evalSymm}}
\end{center}
\end{figure}
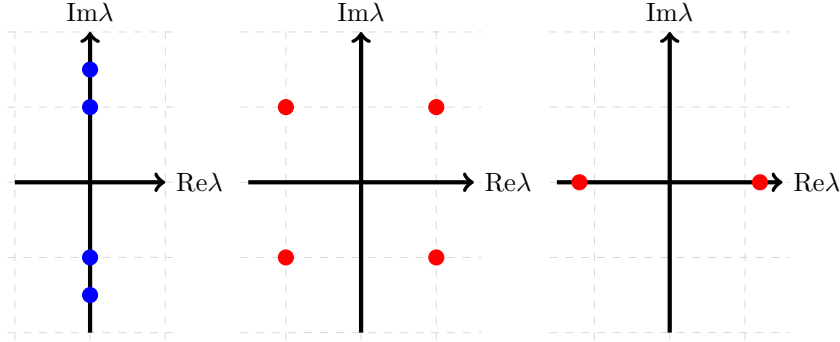

For ease of notation,  we introduce the dispersion relation $\Omega$ in the travelling frame of reference as $\Omega(m+\mu) = \omega(m+\mu) - (m+\mu)V$ with $\lambda(\mu + m) = -i\Omega(m+\mu)$. Since $\lambda$ is purely imaginary when we consider the linear regime, the perturbation will not grow exponentially in time and thus $u^{(0)}(x)$ is spectrally stable. However, as the nonlinearity is increased with increasing $\epsilon$, the eigenvalues which depend continuously on the amplitude of the solution will change and may develop some non-zero real part. Since the equation is Hamiltonian, they will do so symmetrically in the complex plane to conserve the energy, keeping the solution real. The possible configurations of the symmetries in eigenvalues are shown in Figure \ref{fig:evalSymm}. In order to leave the imaginary axis and develop instability, the eigenvalues first have to collide in order to maintain the symmetry of the equation. In Figure \ref{fig:evalSymm}, even if eigenvalues move and collide, they do not necessarily leave the imaginary axis as shown in the left panel. This implies a necessary condition for instability is collisions of eigenvalues for different modes $m$ and $n$ in a perturbation given by
\begin{align}
\lambda(\mu + m) = \lambda(\mu + n).
\label{eq:colGen}
\end{align}

Also in the linear regime (considering the $O(\delta)$ term when substituting \eqref{eq:deltaExp} into \eqref{eq:genForm}  with  $|u^{(0)}(x)| \rightarrow 0$), we can explicitly write the Hamiltonian of the system as
\begin{align}
H_{\text{lin}} =\int_0^{L} \frac{1}{2}\left(\sum_{n=1}^N (-1)^{n} C_{2n+1} (u^{(1)}_{nx})^2 + V (u^{(1)})^2\right) dx,
\label{eq:Hlin}
\end{align}
with 
\begin{align}
0 = \partial_x \frac{\delta H_{\text{lin}}}{\delta u^{(1)}}.
\end{align}
An unstable solution has to conserve energy given by \eqref{eq:Hlin}. This implies that for a collision of eigenvalues arising from two different modes, for every mode that is contributing positively to the Hamiltonian, there needs to be a negatively contributing mode as well. This contribution of eigenvalues to the Hamiltonian (known as their \textbf{signature}) is simply given by the sign of the Hamiltonian. The signature is derived from \eqref{eq:Hlin} by substituting $u^{(1)} \sim e^{i(\mu + m)x}$ to obtain
\begin{align}
\text{sign}(H_{\text{lin}}) = \text{sign}\left(\sum_{m=1}^N(-1)^mC_{2m+1}(i(\mu + m))^{2m} + V \right).
\end{align}
Using the definition of the dispersion relation in the moving frame and dividing by $i$, we can write the sign of the Hamiltonian as 
\begin{align}
\text{sign}(H_{\text{lin}}) = \text{sign}\left(\frac{\Omega(\mu + m)}{\mu + m}\right),
\end{align}
where we have used \eqref{eq:dispersion} and the definition of the dispersion relation incorporating the travelling frame of reference. With more algebra described in \cite{DT15, KDT18}, we can introduce $s$ which will govern if two colliding eigenvalues for modes $m$ and $n$ will have opposing signature as
\begin{align}
s = (\mu +m)(\mu+n) < 0.
\label{eq:sigCond}
\end{align}

To reduce the number of unknowns in \eqref{eq:sigCond}, we set $(\mu+m) \rightarrow \mu$ therefore letting $n \rightarrow (n-m)$, shifting the focus instead on the difference in Fourier modes of the perturbation. This implies that if we wish to consider when the periodic travelling wave solutions are unstable to perturbations of the form shown in \eqref{eq:pertForm}, then we need examine the \textbf{collision condition} 
\begin{align}
\lambda(\mu)=\lambda(\mu+n)
\label{eq:collision}
\end{align}
as well as the corresponding combination of \textbf{signatures of colliding eigenvalues} given by
\begin{align}
s = \mu(\mu+n)
\label{eq:signature}
\end{align}
In the following sections, we show this can be further simplified to one condition using a reduced order polynomial of degree $N$ and examine where the polynomial has real roots  thereby meeting the necessary conditions for instability.

\section{General Methodology}\label{sec:methods}
In general, if we are given a polynomial with $p(\mu)=\mu^N$ with $N$ odd (for example one term in a dispersion relation), then a collision of eigenvalues is of the form
\begin{align}
p(\mu+n)-p(\mu) = 0.
\label{eq:genPolyCol}
\end{align}
Setting $s = \mu(\mu+n)$, we can equivalently write the collision condition as a reduced-order polynomial $q(s,n)$ of order $\frac{N-1}{2}$ that is indirectly dependent on the Floquet parameter $\mu$ as
\begin{align}
q(s,n) = \sum_{i=0}^{\frac{N-1}{2}}a_{i,N-2i}s^in^{N-2i}.
\label{eq:sGen}
\end{align}
The coefficients can be computed recursively as
\begin{align}
a_{i,j} = 
\begin{cases}
{N \choose j} \  &\text{for} \ i=0,j=2,...,N, \\
a_{i-1,j+1} - a_{i,j+1} \ &\text{for} \ i=1,...,\frac{N-1}{2}, j=1,...,N-2i,\\
0 \ &\text{otherwise}.
\end{cases}
\end{align}
Rewriting the collision condition as a signature condition is always possible as shown by Kollar et al. in \cite{KDT18}. In the following section we will focus on the simplicity of constructing this polynomial for the signature. The main consequence of being able to rewrite the polynomial of lower order, is that it simplifies the equation and the number of roots we have to consider. From \eqref{eq:signature}, we can solve for the Floquet parameter as 
\begin{align}
\mu = \frac{1}{2}\left(-n \pm \sqrt{n^2+4 s} \right).
\label{eq:muSoln}
\end{align}
To satisfy both the collision condition and signature condition for instability while maintaining that perturbations are bounded in space, we need the roots of \eqref{eq:muSoln} to be real and for the signatures to remain opposite, i.e.
\begin{align}
-\frac{n^2}{4} < s < 0.
\label{eq:interval}
\end{align}

Checking that the roots of a polynomial are within a certain interval $I$, in this case given by \eqref{eq:interval}, becomes a relatively straightforward procedure and is in some respect easier than computing exact roots. This can be done using Sturm's theory \cite{SturmBook, T41} via a sequence of polynomials (sometimes known as a Sturm chain). Given a polynomial $g(x) = g_0(x)$ of degree $N$ with real coefficients, a sequence of polynomials of decreasing order is constructed by using the following criteria
\begin{align}
g_1(x) & = \frac{\partial}{\partial x} g_0(x)\ \text{and} \\
g_n(x) & = -\left( g_{n-2}(x)-g_{n-1}(x)\frac{g_{n-2}(x)}{g_{n-1}(x)} \right) = -\text{Rem}(g_{n-2}(x),g_{n-1}(x))
\label{eq:sturmSeq}
\end{align}
where $\frac{g_{n-2}(x)}{g_{n-1}(x)}$ is a polynomial quotient and Rem($g_{n-2}(x),g_{n-1}(x)$) is the remainder. The sequence terminates at $n=N$ when the last term is a constant and therefore independent of $x$. If we are interested in how many real roots $r_n$ occur in the interval $I = (a_i,a_f)$, where $a_i$ and $a_f$ are not themselves roots, then we need to examine the difference in the number of sign changes of the polynomials evaluated at the endpoints of the interval (as shown in \eqref{eq:interval}, in this case $a_i = -n^2/4$ and $a_f = 0$). To obtain the number of real roots in the interval, we subtract the number of sign changes at $a_f$ from the number of sign changes at $a_i$.

To summarise, in order to analyse spectral stability of periodic travelling waves of \eqref{eq:travGenForm} to high-frequency instabilities of the form given by \eqref{eq:pertForm}, we must 
\begin{enumerate}
\item Write the dispersion relation $\omega$ given by the general form in \eqref{eq:dispersion}.
\item Compute the travelling wave speed $V$ for a non-trivial solution.
\item Solve for the polynomial that governs the collision condition of the form \eqref{eq:collision}.
\item Reduce the order of the polynomial by substituting $s=\mu(\mu+n)$.
\item Generate the Sturm sequence of polynomials using \eqref{eq:sturmSeq}.
\item Compute the number of roots in $I$ by examining the number of sign changes in the Sturm sequence of polynomials at each end point and noting the difference.
\end{enumerate}
If the result is that we have no real roots contained in $I$, then the periodic travelling waves are spectrally stable to high-frequency perturbations. In order to show how this method works, we proceed with an example.

\section{Example}\label{sec:example}
In this section we examine an equation of the form
\begin{align}
u_t + \alpha  u_{3x} + \beta u_{5x} + \gamma u_{7x} + f(u)_x = 0,
\label{eq:genFifth}
\end{align}
where $\alpha$, $\beta$ and $\gamma$ are real coefficients and the subscripts represent the number of derivatives of $u(x,t)$ and go through the process outlined in Section 3 to compute the regions of stability, referring to step number in parentheses. In this section, we will keep these as variables however in practice, they are defined by the scaling in the partial differential equation that is being considered. We begin by introducing a travelling frame of reference, moving with speed $V$ and considering a steady-state solution
\begin{align}
\alpha u_{3x} + \beta u_{5x} + \gamma u_{7x} + f(u)_x + V u_x = 0.
\end{align}
The dispersion relation (step 1 in the process) of this equation is given by
\begin{align}
\omega = -\alpha k^3 + \beta k^5 - \gamma k^7.
\label{eq:partDispersion}
\end{align}
Linearizing about a small amplitude solution with $u^{(0)}=\epsilon e^{ikx}$ (where $f(u^{(0)}_x) \approx 0$), we obtain
\begin{align}
 \alpha (ik)^3 + \beta (ik)^5 + \gamma (ik)^7 + V(ik) = 0,
\end{align}
or 
\begin{align}
- \alpha k^2 + \beta k^4 - \gamma k^6 + V = 0.
\end{align}
If we assume the solution we are linearising about is $2 \pi$ periodic, we can show it is symmetric and without loss of generality we can set $k=1$. This gives $V_0 = \alpha - \beta + \gamma$ (completing step 2) as a bifurcation point from which we can compute non-trivial solutions $u^{(0)}(x)$ travelling at speed $V_0$.  We will sub in for $V=V_0$ in the equations from now on. 

The polynomial in terms of $(\mu,n)$ (step 3) for the collision condition is given by
\begin{align}
p(\mu,n) = \gamma (\mu+n)^7 - \beta (\mu+n)^5 + \alpha (\mu+n)^3  - \gamma \mu^9  + \beta \mu^5 -\alpha \mu^3 - (\alpha - \beta + \gamma)n.
\end{align}
The above can be simplified if we set $s = \mu(\mu + n)$. In order to do this, we first note that we can use binomial theorem gives us the polynomial expansion 
\begin{align}
(\mu+n)^N = \sum _{k=0}^{N}{N \choose k}\mu^{N-k}n^{k}.
\end{align}
\begin{table}
\begin{center}
\begin{tabular}{>{$ }l<{$}*{21}{c}}
&&&&&&&&&&1&&&&&&&&\\
&&&&&&&&&1&&1&&&&&&&&\\
&&&&&&&&1&&2&&1&&&&&&&\\
&&&&&&&1&&3&&3&&1&&&&&&\\
&&&&&&1&&4&&6&&4&&1&&&&&\\
&&&&&1&&5&&10&&10&&5&&1&&&&\\
&&&&1&&6&&15&&20&&15&&6&&1&&&\\
&&&1&&7&&21&&35&&35&&21&&7&&1&&\\
\end{tabular}
\caption{Coefficients from the binomial theorem in a Pascal's triangle \label{tab:pascal}.}
\end{center}
\end{table}
The coefficients from the binomial theorem can be computed via Pascal's triangle where each row represents coefficients in a polynomial of degree $N=0,\cdots,7$ shown in Table \ref{tab:pascal} and obtain the collision condition as
\begin{align}
p(\mu,n) = & \gamma (7\mu^6 n+ 21\mu^5 n^2 + 35\mu^4 n^3+35\mu^3 n^4 + 21\mu^2 n^5 + 7\mu n^6+n^7) \nonumber \\
- & \beta \left( 5\mu^4 n + 10 \mu^3 n^2 + 10\mu^2 n^3 + 5 \mu n^4 + n^5 \right) + \alpha ( 3\mu^2n + 3\mu n^2 + n^3 )\nonumber \\   & - \left( \alpha - \beta + \gamma \right()n = 0.
\label{eq:muExpanded}
\end{align}
Just as Pascal's triangle provides an easy way to compute the coefficients of $(\mu+n)^N$ in \eqref{eq:genPolyCol}, we can use a triangular construction to find the coefficients of $q(s,n)$ in \eqref{eq:sGen}. To begin, create a table whose $N$ columns are the coefficients of  $(\mu+n)^N-\mu^N$ beginning with the coefficient of $n^{N}$ and ending with the coefficient of $n^1$.
Row 2 begins with a zero one place to the left of the first column in row 1. Subsequent elements in row 2 are found by computing the difference between row 1 and row 2 in the previous column.
This procedure is repeated until the final row which will have just two elements.
The coefficients in the reduced polynomial for the signature (that is the polynomial which depends on $s=\mu(\mu+n)$) are the first non-zero values in each row (circled in Tables \ref{tab:s7}-\ref{tab:s3} below).  They are given in increasing order of $s$ as labelled in the right-most row.  That is, row 1 gives the coefficient of $s^0n^N$ and row $(N+1)/2$ gives the coefficient of $s^{(N-1)/2}n^1$. Tables \ref{tab:s7} - \ref{tab:s3} show this process explicitly for $N=7,5,3$ respectively.

\begin{table}[h!]
\centering
\begin{tabular}{|cccccccccccccc|c|}
\hline
 \bf{$n^7$} && \bf{$n^6$} && \bf{$n^5$} && \bf{$n^4$} && \bf{$n^3$} && \bf{$n^2$} && \bf{$n^1$} &&  \\ \hline
\circled{1}&&7&&21&&35&&35&&21 &&7 && $s^0$ \\
&&\textcolor{blue}{$\downarrow$}&&\textcolor{blue}{$\downarrow$}&&\textcolor{blue}{$\downarrow$}&&\textcolor{blue}{$\downarrow$}&&\textcolor{blue}{$\downarrow$} &&  && \\
&&0&$\rightarrow$&\circled{7}&$\rightarrow$&14&$\rightarrow$&21&$\rightarrow$&14 &$\rightarrow$&7 && $s^1$ \\
&& && &&\textcolor{blue}{$\downarrow$} &&\textcolor{blue}{$\downarrow$}&&\textcolor{blue}{$\downarrow$} && && \\
&&& &&&0&$\rightarrow$&\circled{14}&$\rightarrow$&7 &$\rightarrow$&7 && $s^2$\\
&& && &&  && &&\textcolor{blue}{$\downarrow$} && && \\
&&& &&& & & & &0 &$\rightarrow$&\circled{7} && $s^3$ \\
\hline
\end{tabular}
\caption{Tabular computation of $(\mu+n)^7-\mu^7=n^7+7sn^5+14s^2n^3+7s^3n$. The coefficients of the reduced polynomial in terms of $s=\mu(\mu+n)$ are given by the circled terms. Downward arrows (in blue) indicate subtraction and arrows to the right (in black) indicate the result of the subtraction. \label{tab:s7}}
\end{table}

\begin{table}[h!]
\centering
\begin{tabular}{|cccccccccc|c|}
\hline
\bf{$n^5$} && \bf{$n^4$} && \bf{$n^3$} && \bf{$n^2$} && \bf{$n^1$} && \\\hline
\circled{1}& &5& &10& &10 & &5 && $s^0$ \\
 &&\textcolor{blue}{$\downarrow$} &&\textcolor{blue}{$\downarrow$}&&\textcolor{blue}{$\downarrow$} && && \\
&&0&$\rightarrow$&\circled{5}&$\rightarrow$&5 &$\rightarrow$&5 && $s^1$\\
 &&  && &&\textcolor{blue}{$\downarrow$} && && \\
&& & & & &0 &$\rightarrow$&\circled{5} && $s^2$\\
\hline
\end{tabular}
\caption{Tabular computation of $(\mu+n)^5-\mu^5=n^5+5sn^3+5s^2n$. The coefficients of the reduced polynomial in terms of $s=\mu(\mu+n)$ are given by the circled terms. Downward arrows (in blue) indicate subtraction and arrows to the right (in black) indicate the result of the subtraction.\label{tab:s5}}
\end{table}

\begin{table}[h!]
\centering
\begin{tabular}{|cccccc|c|}
\hline
\bf{$n^3$} && \bf{$n^2$} && \bf{$n^1$} && \\\hline
\circled{1}&$\rightarrow$&3 &$\rightarrow$&3 && $s^0$\\
 &&\textcolor{blue}{$\downarrow$} && && \\
 & &0 &$\rightarrow$&\circled{3} && $s^1$ \\
\hline
\end{tabular}
\caption{Tabular computation of $(\mu+n)^3-\mu^3=n^3+3sn$. The coefficients of the reduced polynomial in terms of $s=\mu(\mu+n)$ are given by the circled terms. Downward arrows (in blue) indicate subtraction and arrows to the right (in black) indicate the result of the subtraction.\label{tab:s3}}
\end{table}

Finally, combining the results from the Tables \ref{tab:s7} - \ref{tab:s3}, the polynomial for the signature condition (step 4 in the process) is 
\begin{align}
q(s,n) = & -\gamma(n^6+7n^4s+14n^2s^2+7s^3)+\beta(n^4+5n^2s+5s^2) \nonumber \\
&-\alpha(n^2+3s)+\alpha-\beta+\gamma
\label{eq:q753}
\end{align}

We can analyse the roots of \eqref{eq:q753} using Sturm's theory by constructing a sequence of polynomials (this is step 5) in $s$ of the form in \eqref{eq:sturmSeq} as
\begin{align}
p_1(s) = & -\gamma (7 n^4 + 28n^2s + 21s^2) + \beta(5n^2 + 10s)-3\alpha  \label{eq:p1}\\
p_2(s) = & -\frac{2s}{63 \gamma}\left(49 \gamma^2 n^4 - 35\beta\gamma n^2-63 \alpha\gamma+25\beta^2\right) \nonumber \\
& -\frac{1}{63\gamma}\left(35\gamma^2n^6-42\beta\gamma n^4-21\alpha\gamma n^2+25\beta^2 n^2-15\alpha\beta+63\alpha\gamma-63\beta\gamma+63\gamma^2\right) \label{eq:p2}\\
p_3(s) = & -\frac{63\gamma}{4(49\gamma^2 n^4- 35\beta\gamma n^2-63\alpha\gamma+25\beta^2)^2} \left(49\gamma^4n^{12}-196\beta\gamma^3n^{10} - 98\alpha\gamma^3n^8 \right. \nonumber \\ 
& \left. + 322\beta^2\gamma^2n^8 - 126\alpha\beta\gamma^2 n^6 + 1274\alpha\gamma^3n^6-200\beta^3\gamma n^6 - 1274\beta\gamma^3 n^6 + 1274\gamma^4n^6 \right. \nonumber \\
& \left. + 441\alpha^2\gamma^2 n^4 - 210\alpha\beta^2 \gamma n^4 - 1176\alpha\beta\gamma^2 n^4 + 125\beta^4n^4 + 1176\beta^2\gamma^2 n^4 \right. \nonumber \\
&\left. - 1176\beta\gamma^3 n^4 + 630\alpha^2 \beta\gamma n^2 - 2646\alpha^2\gamma^2n^2- 250\alpha\beta^3n^2 + 1050\alpha\beta^2\gamma n^2 \right. \nonumber \\
& \left. + 2646\alpha\beta\gamma^2n^2 - 2646\alpha\gamma^3n^2 - 1050\beta^3\gamma n^2 + 1050\beta^2\gamma^2 n^2-756\alpha^3\gamma \right. \nonumber \\
& \left. + 225\alpha^2\beta^2 +1890\alpha^2\beta\gamma - 1323\alpha^2\gamma^2 -500\alpha\beta^3 -1890\alpha\beta^2\gamma+ 4536\alpha\beta\gamma^2 \right. \nonumber \\
&\left. -2646\alpha\gamma^3 +500\beta^4 -500\beta^3\gamma -1323\beta^2\gamma^2 + 2646\beta\gamma^3-1323\gamma^4 \right)
\label{eq:p3}
\end{align}
Despite the length of the expressions in \eqref{eq:q753}-\eqref{eq:p3}, their sign changes are easy to evaluate for particular $\alpha, \beta, \gamma$ and $s \in (-1/4,0)$. For ease, Table \ref{tab:g0ns} shows the sign changes for $\alpha = 1$, $\beta = 1/4$ and $\gamma = 0$ for $n=1,2,3,4$ which are in complete agreement with results in \cite{TDK18} (this is the final step in the process, step 6). They imply that the perturbations with $n\ge4$ are stable as is the perturbation for $n=1$ since there are no real roots. Note that in cases where $p_j(s) = 0$, we must consider the limit as $s$ approaches the value $0$ or $-n^2/4$ from the correct side to match with the condition in \eqref{eq:interval}.

\begin{table}[h!]
\begin{center}
\begin{tabular}{|l|c|c|}
\hline
$n=1$ & sign($p_j(\!-\!n\!^2\!/\!4)$) & sign($\!p_j(0)$) \\
\hline
$p_0(s)$ & + & +\\
$p_1(s)$ & - & -\\
$p_2(s)$ & + & + \\
\hline
sign &  &  \\
changes & 2 & 2\\
\hline
\end{tabular}
%
\begin{tabular}{|l|c|c|}
\hline
$n=2$ & sign($p_j(\!-\!n\!^2\!/\!4)$) & sign($\!p_j(0)$) \\
\hline
$p_0(s)$ & - & + \\
$p_1(s)$ & - & + \\
$p_2(s)$ & + & + \\
\hline
sign &  &  \\
changes & \! 1 &  \! 0 \\
\hline
\end{tabular}
%
\begin{tabular}{|l|c|c|}
\hline
$n=3$ & sign($p_j(\!-\!n\!^2\!/\!4)$) & sign($\!p_j(0)$) \\
\hline
$p_0(s)$ & - & + \\
$p_1(s)$ & + & + \\
$p_2(s)$ & + & + \\
\hline
sign &  &  \\
changes & 2 & 0 \\
\hline
\end{tabular}
%
\begin{tabular}{|l|c|c|}
\hline
$n=4$ & sign($p_j(\!-\!n\!^2\!/\!4)$) & sign($\!p_j(0)$) \\
\hline
$p_0(s)$ & + & + \\
$p_1(s)$ & + & + \\
$p_2(s)$ & + & + \\
\hline
sign &  &  \\
changes & 0 & 0 \\
\hline
\end{tabular}
\end{center}
\caption{The stability results with $\alpha = 1$, $\beta = 1/4$ and $\gamma = 0$ and $n=1,2,3,4$ (note this is a singular case of \eqref{eq:genFifth}). By subtracting the number of sign changes at $s=-n^2/4$ from the ones at $s=0$ (subtract the total in column 3 from column 2) for each $n$, we get the number of real roots in that interval. Instability is possible for $n=2$ and $n=3$ since there are roots for which $s \in I=(-n^2/4,0)$. We can also conclude the equation is stable to perturbations with $n=1$ and $n\ge4$.\label{tab:g0ns}}
\end{table}

\section{Stability Results}\label{sec:results}
\begin{figure}
\begin{center}
\includegraphics[width=0.49\textwidth]{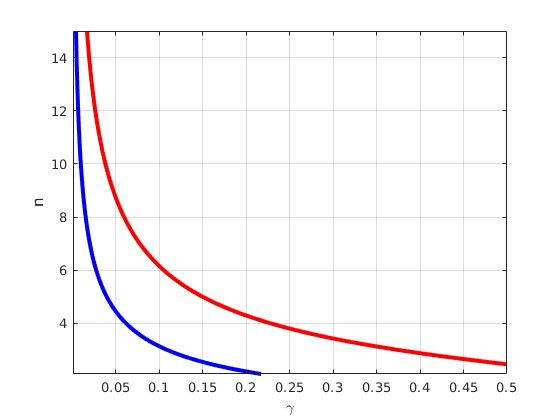} 
\includegraphics[width=0.49\textwidth]{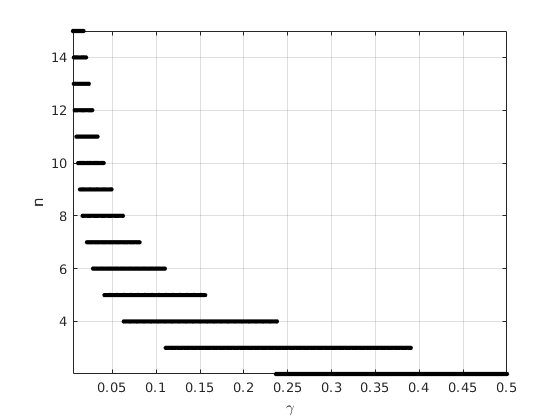}
\caption{Reduction to the two-dimensional system with $\alpha = 0$, giving the instability results for $u_t=Vu_x + \beta u_{5x} + \gamma u_{7x} + \text{nonlinearity}$. \label{fig:a0}}
\end{center}
\end{figure}

\begin{figure}
\includegraphics[width=0.49\textwidth]{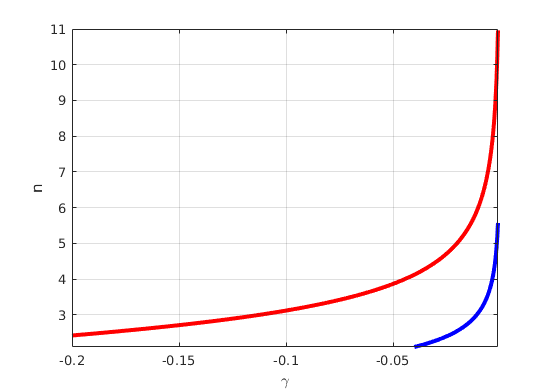} 
\includegraphics[width=0.49\textwidth]{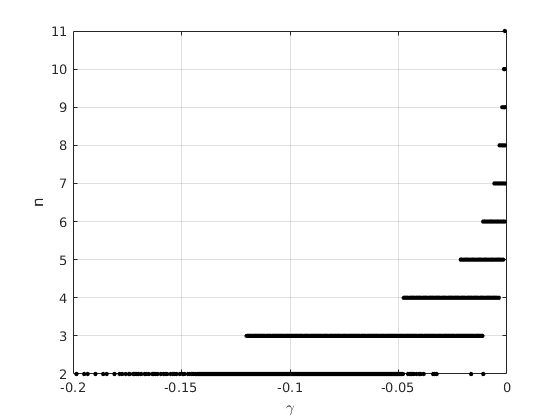}
\caption{Reduction to the two-dimensional system with $\beta = 0$, giving the instability results for $u_t=Vu_x + \alpha  u_{3x} + \gamma u_{7x} + \text{nonlinearity}$.\label{fig:b0}}
\end{figure}

\begin{figure}
\begin{center}
\includegraphics[width=0.49\textwidth]{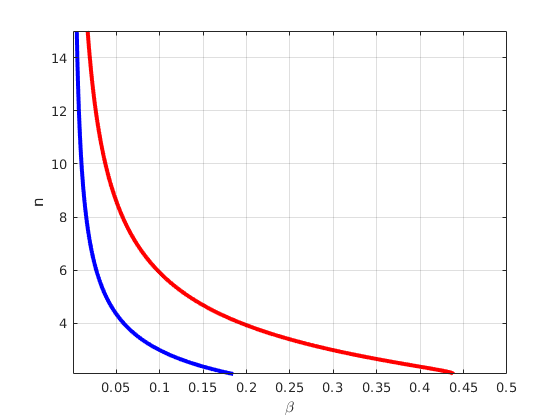} 
\includegraphics[width=0.49\textwidth]{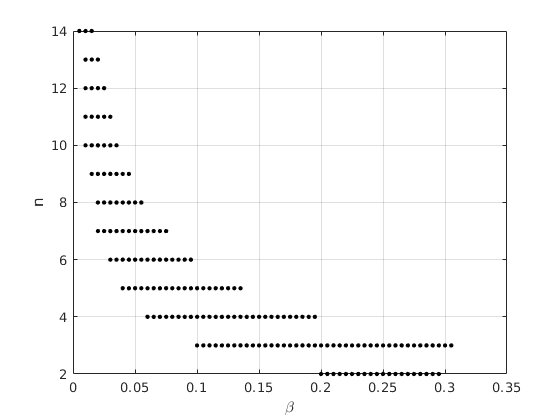}
\caption{Reduction to the two-dimensional system with the singular limit when $\gamma = 0$, giving the instability results for $u_t=Vu_x + \alpha u_{3x} + \beta u_{5x} + \text{nonlinearity}$.\label{fig:g0}}
\end{center}
\end{figure}

Figures \ref{fig:a0} - \ref{fig:g0} show in more detail the stable and unstable regions in two-dimensions for PDEs with only one free parameter (setting one of the parameters in the PDE to zero). In Figure \ref{fig:a0}, $\alpha = 0$, $\beta = 1$ and $\gamma$ is a free parameter. The region bounded below by the blue line and above by the red line is where we can have instability and outside of these curves is where the small amplitude solutions are stable with respect to the instabilities considered in this work. In the plot on the right in Figure \ref{fig:a0}, the dots show the unstable regime for integer values of $n$ where \eqref{eq:q753} has roots in the interval $(-n^2/4,0)$. We see that as $\gamma$ decreases, the instabilities occur for larger $n$, indicating the difference in Fourier modes of colliding eigenvalues. Figure \ref{fig:b0} gives the stability regions for $\alpha = 1$, $\beta = 0$ and $\gamma$ as a free parameter. In this case, only $\gamma < 0$ leads to instabilities, but the pattern is similar to the previous figure. Figure \ref{fig:g0} gives the results previously computed in \cite{TDK18} where once again with decreasing $\beta$, the instabilities have an increasing $n$.

\begin{figure}
\begin{center}
\includegraphics[width=0.8\textwidth]{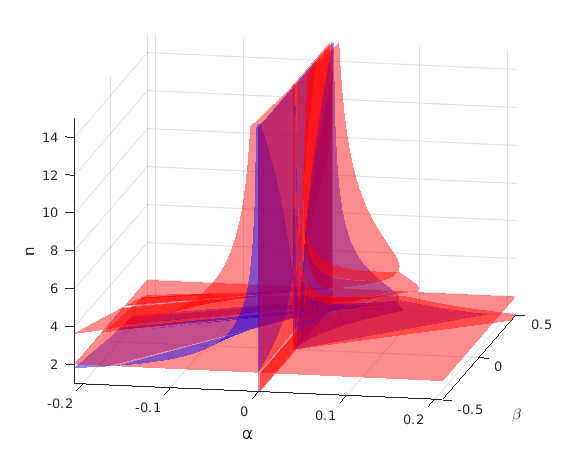} \\
\includegraphics[width=0.8\textwidth]{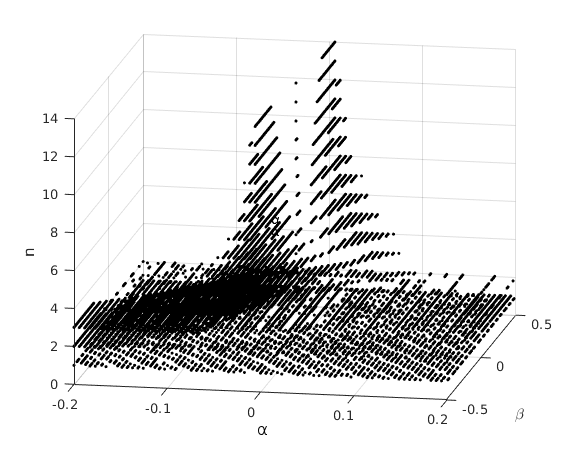}
\caption{On the top, the regions bounded by red and blue curves are those where instabilities can arise. On the bottom, the dots represent possible unstable regions for discrete values of $n$. For both figures, the equations were normalised such that $\gamma = 1$. \label{fig:summary}}
\end{center}
\end{figure}

Figure \ref{fig:summary} summarises the full stability results for the general PDE \eqref{eq:genForm} with $\gamma=1$, which is simply a rescaling of the full equation and does not reduce the degrees of freedom. The regions  between the blue and red lines are possible regions of instability. For clarity, points in the lower plot of Figure \ref{fig:summary} show possible regions of instability and the white space gives the regimes for spectrally stable periodic travelling wave solutions to \eqref{eq:genForm}. This plot shows that most of the regimes of \eqref{eq:genForm} are stable.

\section{Conclusion}\label{sec:conclusion}
In this work, we describe a systematic way to fully characterise spectral stability regions of travelling wave solutions of a dispersive, Hamiltonian PDE subject to high-frequency instabilities. This method shows explicitly how two necessary conditions can be merged into one and a systematic way to analyse the reality of its roots. It relies on reducing the polynomial derived from the dispersion relation describing collisions of eigenvalues of degree $2N+1$, to a polynomial for the signature condition of degree $N$. This polynomial can be constructed using a triangle of coefficients as is illustrated using an example of a PDE containing three linear dispersive terms with general coefficients. If this reduced-order polynomial has roots in a given interval $I = (-n^2/4,0)$, which can be determined using Sturm's theory, then the necessary criteria for an instability is met. This methodology can be used on any dispersive, Hamiltonian partial differential equation. Sturm's theory has also been implemented in Maple and can be accessed through the commands \texttt{sturm} and \texttt{sturmseq}.

There are two drawbacks to this method. One is that it can only be used if the sign of the Hamiltonian is definite, hence the restriction to high-frequency instabilities is made. It also relies on the underlying equations having a Hamiltonian and hence a four-fold symmetry in the complex eigenvalue plane. Since many physical systems are Hamiltonian, there is a large number of applications of this method (for more examples, see \cite{DT15}), which also includes Euler equations describing water waves. 

\section{Acknowledgements}
We thank B. Deconinck, R. Koll{\'a}r and D. Ambrose for insightful discussions. We wish to thank Casa Mathem{\'a}tica Oaxaca,  the Erwin Schr\"{o}dinger Institute and Institute for Computational and Experimental Research in Mathematics (ICERM) for their hospitality during the development of the ideas for this work.

\bibliographystyle{plain}

\bibliography{mybib}

\end{document}